\documentclass[a4paper,11pt]{amsart}

\usepackage{amsmath,amsthm,amssymb}
\usepackage[mathscr]{eucal}
\usepackage{cite}
\usepackage[bookmarks,bookmarksnumbered]{hyperref}

\usepackage[all]{xy}

\numberwithin{equation}{section}

\theoremstyle{plain}
\newtheorem{theorem}{Theorem}[section]

\newtheorem{lemma}[theorem]{Lemma}

\newtheorem{corollary}[theorem]{Corollary}
\theoremstyle{definition}
\newtheorem{definition}[theorem]{Definition}

\begin{document}

\title[Weak containment rigidity]{Weak containment rigidity for distal actions}

\author[Adrian Ioana and Robin Tucker-Drob]{Adrian Ioana and Robin Tucker-Drob}
\thanks{A.I. was partially supported by NSF Grant DMS \#1161047,  NSF Career Grant DMS \#1253402, and a Sloan Foundation Fellowship. R.T.D. was partially supported by NSF Grant DMS \#1303921 and by NSF Grant DMS \#1600904}
\address{A.I. Mathematics Department; University of California, San Diego, CA 90095-1555 (United States)
R.T.D. Department of Mathematics, Mailstop 3368, Texas A{\&}M University, College Station, TX 77843-3368, USA}
\email{aioana@ucsd.edu, rtuckerd@math.tamu.edu}

\begin{abstract}
We prove that if a measure distal action $\alpha$ of a countable group $\Gamma$ is weakly contained in a strongly ergodic probability measure preserving action $\beta$ of $\Gamma$, then $\alpha$ is a factor of $\beta$. In particular, this applies when $\alpha$ is a compact action.

As a consequence, we show that the weak equivalence class of any strongly ergodic action completely remembers the weak isomorphism class of the maximal distal factor arising in the Furstenberg-Zimmer Structure Theorem.
\end{abstract}

\maketitle

\section{Introduction}
The notion of weak containment for group actions was introduced by A. Kechris \cite{Ke10} as an analogue of the notion of weak containment for unitary representations. Let $\Gamma\curvearrowright^{\alpha} (X,\mu)$ and $\Gamma\curvearrowright^{\beta}(Y,\nu)$ be two probability measure preserving  (p.m.p.) actions of a countable group $\Gamma$. Then $\alpha$ is said to be {\it weakly contained} in $\beta$ (in symbols, $\alpha\prec\beta$) if for any finite set $S\subset\Gamma$, finite measurable partition $\{A_i\}_{i=1}^n$ of $X$, and $\varepsilon>0$, we can find a measurable partition $\{B_i\}_{i=1}^n$ of $Y$ such that for all $\gamma\in S$ and $i,j\in\{1,...,n\}$ we have $$|\mu(\gamma A_i\cap A_j)-\nu(\gamma B_i\cap B_j)|<\varepsilon.$$
If $\alpha\prec\beta$ and $\beta\prec\alpha$, we say that $\alpha$ is {\it weakly equivalent} to $\beta$.

We say that $\alpha$ is a {\it factor} of $\beta$, or that $\beta$ is an {\it extension} of $\alpha$, if there exists a measurable, measure preserving map $\theta:Y\rightarrow X$ such that $\theta(\gamma y)=\gamma\theta(y)$, for all $\gamma\in\Gamma$ and almost every $y\in Y$. The map $\theta$ is called a {\it factor map} or an {\it extension}. If in addition there is a conull set $Y_0\subseteq Y$ such that $\theta$ is one-to-one on $Y_0$, then $\theta$ is called an {\it isomorphism} and we say that $\alpha$ is {\it isomorphic} to $\beta$. The actions $\alpha$ and $\beta$ are said to be {\it weakly isomorphic} if each is a factor of the other. 

As the terminology suggests, if $\alpha$ is a factor of $\beta$, then $\alpha$ is weakly contained in $\beta$. The main goal of this note is to establish a rigidity result which provides a general instance when the converse holds.

\begin{theorem}\label{main} Let $\Gamma$ be a countable group, $\Gamma\curvearrowright^{\alpha} (X,\mu)$ be a measure distal p.m.p.\ action, and $\Gamma\curvearrowright^{\beta}(Y,\nu)$ be a strongly ergodic p.m.p.\ action.

If $\alpha$ is weakly contained in $\beta$, then $\alpha$ is a factor of $\beta$.
In particular, if a compact action $\alpha$ is weakly contained in a strongly ergodic action $\beta$, then $\alpha$ is a factor of $\beta$.
\end{theorem}

Before recalling the notions involved in Theorem \ref{main}, let us put it into context and outline its proof.

Weak containment and weak equivalence have received much attention since their introduction. In \cite{Ke10}, A. Kechris shows that cost varies monotonically with weak containment and in \cite{Ke12} Kechris uses this monotonicity to obtain a new proof that free groups have fixed price. Several other measurable combinatorial parameters of actions are known to respect weak containment and hence are invariants of weak equivalence; see \cite{AE11, CK13, CKTD12}.

%

In \cite{AW13}, M. Ab\'{e}rt and B. Weiss exhibit a remarkable anti-rigidity phenomenon for weak containment by showing that every free p.m.p.\ action of $\Gamma$ weakly contains the Bernoulli action over an atomless base space. The Ab\'{e}rt-Weiss Theorem was extended in \cite{TD15} and used to show that every weak equivalence class contains ``unclassifiably many'' isomorphism classes of actions, thus ruling out the possibility of weak equivalence superrigidity. These anti-rigidity results stand in marked contrast to the rigidity exhibited in Theorem \ref{main} and Corollary \ref{maincor} below, and suggest that Theorem \ref{main} and Corollary \ref{maincor} are likely optimal.

The first and, thus far, only rigidity results for weak containment  were obtained by M. Ab\'{e}rt and G. Elek. They prove that if a finite action $\alpha$ (i.e. an action on a finite probability space) is weakly contained in a strongly ergodic action $\beta$, then $\alpha$ is a factor of $\beta$ \cite[Theorem 1]{AE10}. From this, they deduce that if two strongly ergodic profinite actions are weakly equivalent, then they are isomorphic \cite[Theorem 2]{AE10}. Recall that a p.m.p. action is called {\it profinite} if it is an inverse limit of finite actions.  Since any profinite action is compact, Theorem \ref{main} and Corollary \ref{maincor} below recover the results of \cite{AE10}.

Our results are new for compact non-profinite actions, and in particular for translation actions $\Gamma\curvearrowright (K,m_K)$ on connected compact groups $K$. Note that the approach of \cite{AE10} relies on the fact that in the case of profinite actions $\Gamma\curvearrowright (X,\mu)$, there are ``many" $\Gamma$-invariant measurable partitions of $X$. As such, it does not apply to translation actions on connected compact groups, since these actions admit no non-trivial invariant measurable partitions.

Instead, the proof of Theorem \ref{main} relies on a new, elementary approach which we briefly outline in the case $\alpha$ is a compact action. Let $d$ be a metric on $X$ such that $(X,d)$ is a compact metric space on which $\Gamma$ acts isometrically (see the definition of compact actions given below).
Next, assuming that $\alpha\prec\beta$, we find a sequence $\theta_n:Y\rightarrow X$ of almost $\Gamma$-equivariant measurable maps. Since $\beta$ is strongly ergodic,  the maps $y\mapsto d(\theta_m(y),\theta_n(y))$ must be asymptotically constant, as $m,n\rightarrow\infty$. When coupled with the compactness of the metric space $(X,d)$, this forces a subsequence of $\{\theta_n\}$ to converge almost everywhere. The limit map $\theta:Y\rightarrow X$ is $\Gamma$-equivariant, and hence realizes $\alpha$ as a factor of $\beta$.

Recall that a p.m.p. action $\Gamma\curvearrowright (Y,\nu)$ is called {\it strongly ergodic} if any sequence $A_n\subseteq Y$ of measurable sets satisfying $\nu(\gamma A_n\triangle A_n)\rightarrow 0$, for all $\gamma\in\Gamma$, is trivial, i.e. $\nu(A_n)(1-\nu(A_n))\rightarrow 0$.

 An extension $(X,\mu ) \rightarrow (X_0,\mu _0 )$ of ergodic p.m.p. actions $\Gamma\curvearrowright (X,\mu )$ and $\Gamma \curvearrowright (X_0,\mu _0 )$ is called {\it compact} (or {\it isometric}) if it is isomorphic to a homogeneous skew-product extension, i.e., if there exist a compact group $K$, a closed subgroup $L<K$, and a measurable cocycle $w: \Gamma \times X _0 \rightarrow K$, such that $\Gamma \curvearrowright (X,\mu )$ is isomorphic to $\Gamma\curvearrowright (X_0,\mu _0 )\otimes (K/L , m_{K/L})$, where $m_{K/L}$ is the unique $K$-invariant Borel probability measure on $K/L$, the action is given by $\gamma (x,kL) = (\gamma x, w(\gamma ,x)kL )$ ($x\in X_0$, $kL\in K/L$), and where the map $(X,\mu ) \rightarrow (X_0,\mu _0 )$ corresponds to the projection map $(X_0,\mu _0 )\otimes (K/L , m_{K/L}) \rightarrow (X_0,\mu _0)$. 

The ergodic action $\Gamma \curvearrowright (X,\mu )$ is called {\it compact} if the extension $(X,\mu )\rightarrow (\{\bullet \} , \delta _{\bullet})$, over a point, is a compact extension. Equivalently, this means that the action $\Gamma \curvearrowright (X,\mu )$ is isomorphic to an action of the form $\Gamma \curvearrowright (K/L , m_{K/L})$, where $\Gamma$ acts by translation on $K/L$ via a homomorphism $\Gamma \rightarrow K$ with dense image. In particular, in this case $X$ can be endowed with a metric $d$ such that $(X,d)$ is a compact metric space on which $\Gamma$ acts isometrically.

Let $(X,\mu )\xrightarrow{\varphi} (X_0,\mu _0 )$ be an extension of ergodic p.m.p.\ actions $\Gamma\curvearrowright ^{\alpha} (X,\mu )$ and $\Gamma \curvearrowright ^{\alpha _0} (X_0,\mu _0 )$. An {\it intermediate extension of $\alpha _0$ within $\alpha$} is a pair of extensions, $(X,\mu )\xrightarrow{\varphi _1} (X_1,\mu _1)\xrightarrow{\varphi _{0}} (X_0, \mu _0)$, such that $\varphi _{0}\circ \varphi _1 = \varphi$. The {\it largest intermediate compact extension of $\alpha _0$ within $\alpha$} is the (essentially unique) intermediate extension $(X,\mu )\xrightarrow{\varphi _1} (X_1,\mu _1)\xrightarrow{\varphi _{0}} (X_0, \mu _0)$ of $\alpha _0$ within $\alpha$ satisfying
\begin{itemize}
\item[(a)] The extension $(X_1,\mu _1 ) \xrightarrow{\varphi _{0}} (X_0, \mu _0 )$ is compact;
\item[(b)] Given any other intermediate extension $(X, \mu ) \xrightarrow{\psi _1} (Y_1,\nu _1 )\xrightarrow{\psi _{0}} (X_0 , \mu _0)$ with $(Y_1,\nu _1 )\xrightarrow{\psi _{0}} (X_0 , \mu _0)$ compact, there exists a factor map $(X_1, \mu _1) \xrightarrow{\theta} (Y_1 ,\nu _1 )$ such that $\theta \circ \varphi _1 = \psi _1$ and $\psi _{0}\circ \theta = \varphi _{0}$.
\end{itemize}
See, e.g., \cite[Chapter 9]{Gl03} for a detailed treatment.

The {\it distal tower} associated to an ergodic p.m.p.\ action $\Gamma \curvearrowright ^{\alpha} (X,\mu )$ is the directed family $(\Gamma \curvearrowright ^{\alpha _{\zeta}} (X_{\zeta}, \mu _{\zeta}))_{\zeta <\omega _1}$ of factors of $\alpha$, satisfying:
\begin{align}
\label{eqn:distal1} &\text{The action }\alpha _0 \text{ is the trivial action on a point mass}; \\
\label{eqn:distal2} &\text{The action }\alpha _{\zeta +1}\text{ is the largest intermediate compact extension of }\alpha _{\zeta}\text{ within }\alpha ;\\
\label{eqn:distal3} &\text{For limit ordinals }\zeta ,\text{ the action }\alpha _{\zeta}\text{ is the inverse limit of }(\alpha _{\zeta '})_{\zeta '<\zeta} .
\end{align}
The least countable ordinal $\eta$ for which $\alpha _{\eta +1}=\alpha _{\eta}$ is called the {\it order} of the tower. The action $\alpha$ is said to be {\it (measure) distal} if $\alpha = \alpha _{\zeta}$ for some $\zeta <\omega _1$. The Furstenberg-Zimmer Structure Theorem \cite{Zi76b, Fu77} states that every ergodic p.m.p.\ action $\alpha$ of $\Gamma$ has a unique maximal distal factor -- namely $\alpha _{\eta}$, where $\eta$ is the order of the distal tower associated to $\alpha$ -- and that $\alpha$ is relatively weakly mixing over this factor.

In the rest of the introduction we present several consequences of Theorem \ref{main} and its proof. The first is an immediate corollary of Theorem \ref{main}.

\begin{corollary}\label{maincor1}
Let $\Gamma\curvearrowright^{\alpha}(X,\mu)$ and $\Gamma\curvearrowright^{\beta}(Y,\nu)$ be p.m.p. actions of a countable group $\Gamma$.
Assume that $\beta$ is strongly ergodic.

If $\alpha$ is weakly contained in $\beta$, then the maximal distal factor of $\alpha$ is a factor of $\beta$.
\end{corollary}

The second consequence of Theorem \ref{main} states that, for strongly ergodic actions, the weak isomorphism class of the maximal distal factor is an invariant of weak equivalence.

\begin{corollary}\label{maincor}
Let $\Gamma\curvearrowright^{\alpha}(X,\mu)$ and $\Gamma\curvearrowright^{\beta}(Y,\nu)$ be p.m.p. actions of a countable group $\Gamma$.
Assume that $\beta$ is strongly ergodic.

If $\alpha$ is weakly equivalent to $\beta$, then the maximal distal factors of $\alpha$ and $\beta$ are weakly isomorphic.

Moreover, if $\alpha$ and $\beta$ are compact, and $\alpha$ is weakly equivalent to $\beta$, then $\alpha$ is isomorphic to $\beta$. \end{corollary}

{\it Remark.} The moreover part of Corollary \ref{maincor} is a consequence of the first assertion. Indeed, if two compact actions are weakly isomorphic, then they are in fact isomorphic (see Lemma \ref{lem:weakisom}). However, this fact does not extend to general distal actions; see \cite{Le89} for an example of two distal actions of $\mathbb Z$ which are weakly isomorphic but not isomorphic.
\hfill

The main lemma used in the proof of Theorem \ref{main} (Lemma \ref{equiv}) also implies the following rigidity property for cocycles of strongly ergodic actions with compact targets, generalizing the implication (1) $\Rightarrow$ (2) of \cite[Proposition 2.3]{Sc80} as well as \cite[Lemma J]{Io13}.
Recall that the space $Z^1(\alpha , K )$, of all measurable cocycles of a p.m.p. action $\Gamma \curvearrowright ^{\alpha} (Y, \nu)$ with values in a compact metrizable group $K$, is naturally a Polish space (see e.g., \cite[section 24]{Ke10}). 

\begin{theorem}\label{closed}
Let $\Gamma \curvearrowright ^{\alpha} (Y, \nu)$ be a strongly ergodic p.m.p. action of a countable group $\Gamma$, and $K$ be a compact metrizable group. 

Then the cohomology equivalence relation on $Z^1(\alpha ,K )$ is closed, i.e. is a closed subset of $Z^1(\alpha ,K )\times Z^1(\alpha, K)$. In particular, every cohomology class in $Z^1(\alpha , K )$ is closed.
\end{theorem}

Next, we relate Corollary \ref{maincor} with recent spectral gap results for translation actions on connected compact groups.
A p.m.p. action $\Gamma\curvearrowright(X,\mu)$ has {\it spectral gap} if the unitary representation  $\Gamma\curvearrowright L^2(X)\ominus\mathbb C{\bf 1}$ admits no almost invariant vectors. If an action has spectral gap, then it is strongly ergodic.
 Recent works of J. Bourgain and A. Gamburd \cite{BG06,BG11}, and Y. Benoist and N. de Saxc\'{e} \cite{BdS14} established spectral gap for a wide class of translation actions. More precisely, it was shown that if $K$ is a simple connected compact Lie group and $\Gamma<K$ is a countable dense subgroup generated by matrices with algebraic entries, then the translation action $\Gamma\curvearrowright (K,m_K)$ has spectral gap (see \cite[Theorem 1.2]{BdS14}). Moreover, in this case, the translation action $\Gamma\curvearrowright (K/L ,m_{K/L})$ has spectral gap, for any closed subgroup $L<K$.

 This provides a large family of actions to which Corollary \ref{maincor} applies.
 In particular, it allows us to construct new concrete infinite families of weakly incomparable translation actions of free groups (cf. \cite[Theorem 3]{AE10}).
To this end, let  $a,b$ be integers such that $0<|a|<b$ and $\frac{a}{b}\not=\pm\frac{1}{2}$. Put $c=b^2-a^2$ and let $\Gamma$ be the subgroup of $K=SO_3(\mathbb R)$ generated by the following rotations: $$A=\begin{pmatrix} \frac{a}{b}&-\frac{\sqrt{c}}{b}&0\\\frac{\sqrt{c}}{b}&\frac{a}{b}&0\\0&0&1 \end{pmatrix}\;\;\;\text{and}\;\;\; B=\begin{pmatrix} 1&0&0\\0&\frac{a}{b}&-\frac{\sqrt{c}}{b}\\0&\frac{\sqrt{c}}{b}&\frac{a}{b}\end{pmatrix}.$$

Denote by $\alpha_{(a,b)}$ the associated  translation action $\Gamma\curvearrowright^{\alpha_{(a,b)}}(K,m_K)$. For $2\leq n\leq+\infty$, denote by $\Gamma_n$ the group generated by $\{B^kAB^{-k}|0\leq k\leq n-1\}$, and by $\alpha_{(a,b)}^n$ the restriction of $\alpha_{(a,b)}$ to $\Gamma_n$.
Since $\frac{a}{b}\not=\pm\frac{1}{2}$, \cite{Sw94} shows that $\Gamma$ is isomorphic to $\mathbb F_2$, and thus $\Gamma_n$ is isomorphic to $\mathbb F_n$.

\begin{corollary}\label{cor3} Let $(a,b), (a',b')$ be two pairs of integers as above, and $2\leq n\leq+\infty$.

If $\frac{a}{b}\not=\frac{a'}{b'}$, then $\alpha_{(a,b)}\nprec\alpha_{(a',b')}$ and $\alpha_{(a,b)}^n\nprec\alpha_{(a',b')}^n$.

\end{corollary}

L. Bowen recently proved that if $\beta$ is any essentially free p.m.p. action of a free group $\Gamma=\mathbb F_n$, for some $2\leq n\leq+\infty$, then its orbit equivalence class is weakly dense in the space of all p.m.p. actions of $\Gamma$ (see \cite[Theorem 1.1]{Bo13}). In other words, for any p.m.p. action $\alpha$ of $\Gamma$,  there exists a sequence $\{\beta_k\}$ of actions of $\Gamma$ which are orbit equivalent to $\beta$ and converge to $\alpha$, in the weak topology. As a consequence, $\alpha$ is weakly contained in the infinite direct product action $\times_{k=1}^{\infty}\beta_k$.

In view of this result, it is natural to wonder whether the sequence $\{\beta_k\}$ can be taken constant, that is whether $\alpha$ is weakly contained in some action which is orbit equivalent to $\beta$. By combining Theorem \ref{main} with a result of I. Chifan, S. Popa and O. Sizemore \cite{CPS11}, we are able to show that this is not the case.
More generally, we have:

\begin{corollary}\label{cor4}
Let  $\Gamma\curvearrowright^{\alpha}(X,\mu)$ be  an essentially free ergodic compact p.m.p. action of a countable non-amenable group $\Gamma$. Let $\Lambda\curvearrowright^{\beta}(Y,\nu)=(Z,\rho)^{\Lambda}$ be a Bernoulli action of a countable group $\Lambda$.

Then there does not exist a p.m.p. action $\Gamma\curvearrowright^{\sigma} (Y,\nu)$ such that

\begin{itemize}
\item $\alpha$ is weakly contained in $\sigma$, and
\item $\sigma(\Gamma)(y)\subset\beta(\Lambda)(y)$, for almost every $y\in Y$.
\end{itemize}
\end{corollary}

\subsection*{Acknowledgement} We are very grateful to Eli Glasner and Benjy Weiss for pointing out an error in the first version of the paper, and to Alekos Kechris for helpful remarks. We would also like to thank the anonymous referee for several suggestions that helped improve the exposition.

\section{Proofs}

In this section we prove the results stated in the introduction. We start with some terminology related to weak containment. All groups $\Gamma$ and $\Lambda$ considered below are assumed countable.

\begin{definition}
Let $(Y,\nu )$ be a probability space and let $X$ be a measurable space.
\begin{itemize}
\item[(i)] A sequence $\theta _n :Y \rightarrow X$, $n\in \mathbb{N}$, of measurable maps is said to {\it converge weakly} to the measurable map $\theta : Y\rightarrow X$, if $\nu (\theta _n ^{-1}(A)\triangle \theta ^{-1}(A))\rightarrow 0$ for all measurable subsets $A\subseteq X$.

\item[(ii)] Let $\Gamma \curvearrowright ^{\alpha}(Y,\nu )$ and $\Gamma \curvearrowright ^{\beta}X$ be measurable actions of $\Gamma$ with $\alpha$ probability measure preserving. A sequence $\theta _n :Y\rightarrow X$, $n\in \mathbb{N}$, of measurable maps is said be {\it asymptotically equivariant} for the actions if $\nu(\theta_n^{-1}(\gamma ^{-1}A)\triangle \gamma ^{-1}\theta_n^{-1}(A))\rightarrow 0$ for all $\gamma\in\Gamma$ and measurable subsets $A\subseteq X$.
\end{itemize}
\end{definition}

We now have the following useful characterization of weak containment.

\begin{lemma}\emph{\cite{Ke10}}\label{lemma1}
Let $\Gamma\curvearrowright^{\alpha} (X,\mu)$ and $\Gamma\curvearrowright^{\beta}(Y,\nu)$ be p.m.p. actions.

Then $\alpha\prec\beta$ if and only if there is a sequence $\theta_n:(Y,\nu)\rightarrow (X,\mu)$, $n\in \mathbb{N}$, of measure preserving maps which are asymptotically equivariant.
\end{lemma}

{\it Proof.}   See the proof of \cite[Proposition 10.1]{Ke10}. \hfill$\blacksquare$

\begin{lemma}\label{lemma2}
Let $(Y,\nu )$ be a probability space, let $(X,d)$ be a Polish metric space with $d\leq 1$, and let $\mu$ be a Borel probability measure on $X$.
\begin{enumerate}
\item Let $(\psi _n)_{n\in\mathbb{N}}$ and $(\varphi _n)_{n\in\mathbb{N}}$ be two sequences of measure preserving maps from $(Y,\nu )$ to $(X, \mu )$. Then $\int _Y d(\psi _n (y), \varphi _n (y))\, d\nu \rightarrow 0$ if and only if $\nu (\psi _n^{-1}(A)\triangle \varphi _n ^{-1}(A))\rightarrow 0$ for all measurable subsets $A\subseteq X$.

\item Let $\Gamma\curvearrowright^{\alpha}(X,\mu)$ and $\Gamma\curvearrowright^{\beta}(Y,\nu)$ be p.m.p. actions and let $\theta _n : (Y,\nu )\rightarrow (X,\mu )$, $n\in\mathbb{N}$, be a sequence of measure preserving maps. The sequence $(\theta _n )_{n\in\mathbb{N}}$ is asymptotically equivariant if and only if $\int _Y d(\theta _n(\gamma y), \gamma \theta _n (y))\, d\nu \rightarrow 0$ for each $\gamma \in \Gamma$. Moreover, if $(\theta _n)_{n\in \mathbb{N}}$ is asymptotically equivariant then we may find a subsequence $(\theta _{n_k})_{k\in\mathbb{N}}$ such that $d(\theta_{n_k}(\gamma y),\gamma \theta_{n_k}(y))\rightarrow 0$, for all $\gamma \in\Gamma$ and almost every $y\in Y$.
\end{enumerate}
\end{lemma}

{\it Proof.} (1): Assume first that $\nu (\psi _n^{-1}(A)\triangle \varphi _n ^{-1}(A))\rightarrow 0$ for all measurable subsets $A\subseteq X$. Fix $\varepsilon>0$. Let $\{A_i\}_{i=1}^k$ be a collection of pairwise disjoint Borel subsets of $X$ such that each $A_i$ has diameter at most $\varepsilon$ and $\nu(\bigcup _{i=1}^k A_i )> 1- \varepsilon$. Put $A_0=X\setminus \bigcup _{i=1}^k A_i$. If $i\not=j$, then $$\{y\in Y|\psi _n(y)\in A_i, \varphi _n(y)\in A_j\}\subset \psi _n ^{-1} (A_i)\setminus \varphi _n ^{-1}(A_i).$$ Hence $\nu(\{y\in Y|\psi _n (y)\in A_i,\varphi _n(y)\in A_j\})\rightarrow 0$. Thus, we conclude that $$\nu\big(\bigcup_{i=0}^k\{y\in Y|\psi _n (y)\in A_i,\varphi _n(y)\in A_i\}\big)\rightarrow 1 .$$

Since $\nu(\psi _n^{-1}(A_0))=\mu(A_0)\leqslant\varepsilon$, we get $\liminf _n \nu\big(\bigcup_{i=1}^k\{y\in Y|\psi _n (y)\in A_i, \varphi _n (y)\in A_i\}\big)\geq 1-\varepsilon$. On the other hand, if $\psi _n(y),\varphi _n (y)\in A_i$ for some $1\leq i\leq k$, then $d(\psi _n (y) , \varphi _n(y))\leq\varepsilon$. This proves that $\liminf _n \nu(\{y\in Y|d(\psi _n (y), \varphi _n(y))\leq\varepsilon\})\geq 1-\varepsilon$, for every $\varepsilon>0$, and hence $\int _Y d(\psi _n (y), \varphi _n (y))\, d\nu \rightarrow 0$.

Conversely, assume that $\int _Y d(\psi _n(y), \varphi _n(y) )\, d\nu \rightarrow 0$. Fix $A\subseteq X$ measurable, and $\varepsilon >0$. Since $(X,d)$ is Polish, the measure $\mu$ is tight and regular, so we may find $L \subseteq A \subseteq U$ with $L$ compact, $U$ open, and $\mu (U\setminus L) < \varepsilon$. Letting $r=d(L,X\setminus U) > 0$, we have $\psi _n ^{-1}(L)\setminus \varphi _n ^{-1}(U) \subseteq \{ y\in Y | d(\psi _n(y),\varphi _n(y))\geq r \}$, and hence $\nu (\psi _n ^{-1}(L)\setminus \varphi _n ^{-1}(U) ) \rightarrow 0$. Therefore, since $\psi _n$ and $\varphi _n$ are measure preserving, the containment
\[
\psi _n ^{-1}(A)\setminus \varphi _n ^{-1}(A)\subseteq [\psi _n ^{-1}(L)\setminus \varphi _n ^{-1}(U) ] \cup [\psi _n ^{-1}(A\setminus L)] \cup [\varphi _n ^{-1}(U\setminus A)]
\]
implies that $\limsup _n \nu (\psi _n ^{-1}(A)\setminus \varphi _n ^{-1}(A)) \leq \mu (U\setminus L) <\varepsilon$. Since $\varepsilon >0$ was arbitrary we conclude that $\nu (\psi _n ^{-1}(A)\setminus \varphi _n ^{-1}(A)) \rightarrow 0$, and hence $\nu (\psi _n^{-1}(A)\triangle \varphi _n ^{-1}(A))\rightarrow 0$.

The first part of (2) is immediate from (1), and the moreover statement is clear.
\hfill$\blacksquare$

The proof of Theorem \ref{main} relies on the following result.

\begin{lemma}\label{equiv}
Let $\Gamma\curvearrowright ^{\beta} (Y,\nu)$ be a strongly ergodic p.m.p.\ action. Let $(X,d)$ be a compact metric space, $K$ the group of isometries of $X$, and $w:\Gamma\times Y\rightarrow K$ a measurable cocycle. Assume that there exists a sequence $\theta_n:Y\rightarrow X$, $n\in \mathbb{N}$, of measurable maps such that $\int _Y d(\theta_n(\gamma y),w(\gamma ,y)\theta_n(y))\, d\nu (y)\rightarrow 0$, for all $\gamma \in\Gamma$.

Then there exists a measurable map $\theta:Y\rightarrow X$ such that $\theta(\gamma y)=w(\gamma ,y)\theta(y)$, for all $\gamma \in\Gamma$ and almost every $y\in Y$, and there is a subsequence $(\theta _{n_k})_{k\in\mathbb{N}}$ which converges to $\theta$ both weakly and pointwise (almost everywhere).
\end{lemma}

{\it Proof.} We may clearly assume $d\leq 1$. Given two measurable maps $\theta,\theta':Y\rightarrow X$, we define $$\widetilde{d}(\theta,\theta')=\int_{Y}d(\theta(y),\theta'(y))\, d\nu(y).$$

{\bf Claim.} Let $\theta_n:Y\rightarrow X$ be a sequence of measurable maps and assume that for each $\gamma \in \Gamma$ we have $\int _Y d(\theta_n(\gamma y),w(\gamma ,y)\theta_n(y))\, d\nu (y) \rightarrow 0$. Then for every $\varepsilon>0$, there exists $n\geq 1$ such that the set $\{m\geq n|\widetilde{d}(\theta_m,\theta_n)<\varepsilon\}$ is infinite.

Assuming the claim, let us derive the conclusion. By using the claim we can inductively construct a subsequence $\{\theta_{n_k}\}$ of $\{\theta_n\}$ such that $\widetilde{d}(\theta_{n_{k+1}},\theta_{n_k})<\frac{1}{2^k}$, for any $k\geq 1$. This implies that $$\int_{Y}\sum_{k=1}^{\infty}d(\theta_{n_{k+1}}(y),\theta_{n_k}(y))\;\text{d}\nu(y)=\sum_{k=1}^{\infty}\widetilde{d}(\theta_{n_{k+1}},\theta_{n_k})<1.$$

Therefore,  the sequence $\{\theta_{n_k}(y)\}\subset X$ is Cauchy and thus convergent, for almost every $y\in Y$. It is clear that the map $\theta:Y\rightarrow X$ defined as $\theta(y):=\lim_{k\rightarrow\infty}\theta_{n_k}(y)$ satisfies the conclusion.

{\it Proof of the claim.} Suppose that the claim is false. Thus, there exists $\varepsilon_0>0$ such that for all $n\geq 1$, the set $\{m\geq n|\widetilde{d}(\theta_m,\theta_n)<\varepsilon_0\}$ is finite. Then we can inductively find a subsequence $\{\theta_{n_k}\}$ of $\{\theta_n\}$ such that $\widetilde{d}(\theta_{n_k},\theta_{n_l})\geq\varepsilon_0$, for all $l>k\geq 1$.

Since $X$ is compact, it can be covered by finitely many, say $p\geq 1$, balls of radius $\frac{\varepsilon_0}{4}$. Therefore, if $x_0,x_1,...,x_p$ are $p+1$ points in $X$, then $d(x_i,x_j)\leq\frac{\varepsilon_0}{2}$, for some $0\leq i<j\leq p$.

Next, for $m,n\geq 1$, define $f_{m,n}:Y\rightarrow [0,1]$ by letting $f_{m,n}(y)=d(\theta_{m}(y),\theta_{n}(y))$. Then for all $\gamma \in\Gamma$ and $y\in Y$ we have that

\begin{align*}|f_{m,n}(\gamma y)-f_{m,n}(y)|&=|d(\theta_{m}(\gamma y),\theta_{n}(\gamma y))-d(w(\gamma ,y)\theta_{m}(y),w(\gamma ,y)\theta_{n}(y))|\\&\leq d(\theta_m(\gamma y),w(\gamma ,y)\theta_m(y))+d(\theta_n(\gamma y),w(\gamma ,y)\theta_n(y)).\end{align*}

From this it follows that $\int _Y |f_{m,n}(\gamma y)-f_{m,n}(y)| \, d\nu (y) \rightarrow 0$, for all $\gamma \in\Gamma$, as $m,n\rightarrow\infty$. Since $\|f_{m,n}\|_{\infty}\leq 1$, for all $m,n$, the strong ergodicity assumption implies that  $\|f_{m,n}-\int_{Y}f_{m,n}\|_1\rightarrow 0$, as $m,n\rightarrow\infty$. In other words, $$\int_{Y}|d(\theta_{m}(y),\theta_{n}(y))-\widetilde{d}(\theta_m,\theta_n)|\, d\nu(y)\rightarrow 0,\;\;\;\;\text{as $m,n\rightarrow\infty$}.$$

Recalling that $\widetilde{d}(\theta_{n_k},\theta_{n_l})\geq\varepsilon_0$, for all $l>k\geq 1$, it follows that $$\label{ecu}\nu(\{y\in Y|d(\theta_{n_k}(y),\theta_{n_l}(y))>\frac{\varepsilon_0}{2}\})\rightarrow 1,\;\;\;\;\text{as $k,l\rightarrow\infty$ with $l>k$}.$$

This further implies that $$\nu(\{y\in Y|d(\theta_{n_{k+i}}(y),\theta_{n_{k+j}}(y)>\frac{\varepsilon_0}{2},\;\;\text{for all $0\leq i<j\leq p$}\})\rightarrow 1,\;\;\;\;\text{as $k\rightarrow\infty$}.$$
In particular, if $k$ is large enough, then we can find $y\in Y$ such that $d(\theta_{n_{k+i}}(y),\theta_{n_{k+j}}(y))>\frac{\varepsilon_0}{2}$, for all $0\leq i<j\leq p$. This contradicts the choice of $p$,
and finishes the proof of the claim.
\hfill$\blacksquare$

The following consequence of Lemma \ref{equiv}, which might be of independent interest, will be used below to derive Theorem \ref{closed}.

\begin{corollary}\label{approx_hull} Let $\Gamma\curvearrowright (Y,\nu)$ be a strongly ergodic p.m.p. action. Let $K$ a compact metrizable group endowed with a left-right invariant compatible metric $d$. Let $w:\Gamma\times Y\rightarrow K$ be a measurable cocycle, and $L<K$ be a closed subgroup.
Assume that there exists a sequence of measurable maps $\phi_n:Y\rightarrow K$ such that $d(\phi_n(\gamma y)^{-1}w(\gamma ,y)\phi_n(y),L)\rightarrow 0$, for all $\gamma \in\Gamma$ and almost every $y\in Y$.

Then there exists a measurable map $\phi :Y\rightarrow K$ such that $\phi(\gamma y)^{-1}w(\gamma ,y)\phi(y)\in L$, for all $\gamma \in\Gamma$ and almost every $y\in Y$.
\end{corollary}

This result generalizes \cite[Proposition 2.3]{Sc80} which dealt with the case $L=\{e\}$ and $K=\mathbb T$. It also extends \cite[Lemma J]{Io13}, where it was noticed that the proof of \cite{Sc80} can be adapted to more generally treat the case when $L=\{e\}$ and $K$ is any compact metrizable group.

{\it Proof.} Endow $X=K/L$ with the metric $d'(xL,yL):=d(x,yL)=\inf\{d(x,y\ell )|\ell \in L \}.$
Then $(X,d')$ is a compact metric space and the left multiplication action $K\curvearrowright X$ is isometric. Define $\theta_n:Y\rightarrow X$ by letting $\theta_n(y)=\phi_n(y)K$. Then for all $\gamma \in\Gamma$ and almost every $y\in Y$ we have $$\lim\limits_{n\rightarrow\infty}d'(\theta_n(\gamma y),w(\gamma ,y)\theta_n(y))=\lim\limits_{n\rightarrow\infty}d(\phi_n(\gamma y)^{-1}w(\gamma ,y)\phi_n(y),L)=0.$$

By Lemma \ref{equiv} we deduce the existence of a measurable map $\theta:Y\rightarrow X$ such that we have $\theta(\gamma y)=w(\gamma ,y)\theta(y)$, for all $g\in\Gamma$ and almost every $y\in Y$. Let $\pi:X\rightarrow K$ be a Borel map such that $\pi(x)L=x$, for every $x\in X$ (see e.g. \cite[Theorem 12.17]{Ke95}).
Then $\phi:=\pi\circ\theta:Y\rightarrow K$ is a measurable map such that $\theta(y)=\phi(y)L$, for almost every $y\in Y$, and the conclusion follows. \hfill$\blacksquare$

{\it Proof of Theorem \ref{closed}.} Let $d$ be a left-right invariant compatible metric on $K$.
Assume that $(u _n)_{n\in \mathbb{N}}$ and $(v_n)_{n\in \mathbb{N}}$ are sequences in $Z^1 (\alpha , K )$ converging to $u$ and $v$ respectively, and with $u_n$ and $v_n$ cohomologous for each $n\in \mathbb{N}$. We must show that $u$ and $v$ are cohomologous. Let $w, w_n:\Gamma\times Y\rightarrow K\times K$ be the cocycles defined by $$w(\gamma ,y ) = (u(\gamma ,y), v(\gamma , y))\;\;\;\text{and}\;\;\;w_n(\gamma ,y ) = (u_n(\gamma ,y), v_n(\gamma , y)).$$ For each $n$, since $u_n$ and $v_n$ are cohomologous, there exists a measurable map $\theta _n:Y\rightarrow K$ such that $u_n(\gamma, y)=\theta_n(\gamma y)v_n(\gamma, y)\theta_n(y)^{-1}$, for all $\gamma\in\Gamma$ and almost every $y\in Y$. Define $\phi_n:Y\rightarrow K\times K$ by letting $\phi_n(y)=(\theta_n(y),\text{id}_K)$. Let $L=\{(k,k)|k\in K\}$ be the diagonal subgroup of $K\times K$.
Since $\phi_n(\gamma y)^{-1}w_n(\gamma, y)\phi_n(y)=(v_n(\gamma, y), v_n(\gamma, y))\in L$ and $w_n$ converges to $w$, we deduce that $d(\phi_n(\gamma y)^{-1}w(\gamma, y)\phi_n(y), L)\rightarrow 0$, for all $\gamma\in\Gamma$ and almost every $y\in Y$. By Lemma \ref{equiv} we deduce the existence of a map $\phi:Y\rightarrow K\times K$ such that 
$\phi(\gamma y)^{-1}w(\gamma ,y)\phi(y)\in L$, for all $\gamma \in\Gamma$ and almost every $y\in Y$. This clearly implies that $u$ and $v$ are cohomologous.
\hfill$\blacksquare$

\begin{lemma}\label{lem:distal}
Let $\Gamma\curvearrowright ^{\beta} (Y,\nu)$ be a strongly ergodic p.m.p.\ action and let $\Gamma\curvearrowright ^{\alpha} (X,\mu )$ be a distal p.m.p.\ action. Assume that $\alpha \prec \beta$, as witnessed by the asymptotically equivariant sequence $\theta _n :(Y,\nu )\rightarrow (X,\mu )$, $n\in \mathbb{N}$, of measure preserving maps. Then there exists a factor map $\theta : (Y,\nu ) \rightarrow (X,\mu )$ from $\beta$ to $\alpha$, along with a subsequence $(\theta _{n_k})$ which converges weakly to $\theta$.
\end{lemma}

{\it Proof.}  Let $(\Gamma \curvearrowright ^{\alpha _{\zeta}} (X_{\zeta}, \mu _{\zeta}))_{\zeta <\omega _1}$ be the distal tower associated to $\alpha$, as in \eqref{eqn:distal1}-\eqref{eqn:distal3}, and let $\eta <\omega _1$ be the order of the tower. Since $\alpha$ is distal we have $\alpha = \alpha _{\eta}$. We prove the lemma by transfinite induction on $\eta$. If $\eta =0$ then the statement is obvious, so assume that $\eta >0$.

{\bf Case 1: $\eta = \eta _0 +1$ is a successor ordinal.} In this case the extension $(X , \mu )\rightarrow (X_{\eta _0}, \mu _{\eta _0} )$ is compact, so we may assume without loss of generality that $(X,\mu ) = (X_{\eta _0}, \mu _{\eta _0})\otimes (K/L ,m_{K/L})$, and that the action $\alpha$ is of the form $\gamma (x, kL) = (\gamma x , w(\gamma , x)kL )$ for some measurable cocycle $w:\Gamma \times X_{\eta _0} \rightarrow K$. For each $n$ we may write $\theta _n (y) = (\theta _n ^0 (y), \theta _n ^1 (y))$, where $\theta _n ^0 : Y\rightarrow X_{\eta _0}$ and $\theta _n ^1 :Y \rightarrow K/L$ are the compositions of $\theta _n$ with the left and right projections, respectively.  Applying the induction hypothesis to the action $\alpha _{\eta _0}$ and the sequence $(\theta _n ^0 )$, we obtain a subsequence $(\theta _{n_i}^0)$ and a factor map $\theta ^0 :(Y,\nu )\rightarrow (X_{\eta _0},\mu _{\eta _0})$ such that $\theta _{n_i}^0$ converges weakly to $\theta ^0$.

Define $\bar{\theta }_n : (Y,\nu ) \rightarrow (X, \mu )$ by $\bar{\theta }_n(y) =(\theta ^0(y) , \theta _n ^1 (y))$, so that the subsequence $\bar{\theta }_{n_i}$ is asymptotically equivariant. Fix a compatible Polish metric $d_0\leq 1$ on $(X_{\eta _0}, \mu _{\eta _0})$, let $d_{K/L}\leq 1$ be a compatible $K$-invariant metric on $K/L$, and let $d$ be the compatible Polish metric on $X$ given $d((x,kL),(x',k'L))= \tfrac{1}{2}(d_0(x,x')+d_{K/L}(kL,k'L))$. By Lemma \ref{lemma2}, for each $\gamma \in \Gamma$ we have $\lim \limits_{i\rightarrow \infty} \int _Y d (\bar{\theta }_{n_i}(\gamma y), \gamma \bar{\theta }_{n_i}(y))\, d\nu (y) = 0$. Since $\theta ^0$ is $\Gamma$-equivariant this means that for each $\gamma \in \Gamma$ we have
\[
\lim _{i\rightarrow \infty} \int _Y d_{K/L}(\theta ^1_{n_i}(\gamma y), w(\gamma ,\theta ^0(y))\theta ^1_{n_i}(y)) \, d\nu (y) = 0 .
\]
Applying Lemma \ref{equiv}, we obtain a measurable map $\theta ^1 : Y\rightarrow K/L$ with $\theta ^1(\gamma y) = w(\gamma ,\theta ^0(y))\theta ^1(y)$, along with a subsequence $(\theta _{n_i'}^1)$ of $(\theta _{n_i}^1)$ with $\int _Y d_{K/L}(\theta _{n_i'}^1 (y), \theta ^1 (y))\, d\nu \rightarrow 0$. Therefore, the map $\theta :(Y,\nu )\rightarrow (X,\mu )$ defined by $\theta (y)=(\theta ^0(y), \theta ^1(y))$ is a factor map from $\beta$ to $\alpha$. Moreover, by Lemma \ref{lemma2}, the subsequence $(\theta _{n_i'})$ converges weakly to $\theta$.

{\bf Case 2: $\eta$ is a limit ordinal.} Fix a sequence $(\zeta _j)_{j\in \mathbb{N}}$ of ordinals which strictly increase to $\eta$. For each $j<k\in \mathbb{N}$ let $\varphi _{j,k} : (X_{\zeta _k},\mu _{\zeta _k} ) \rightarrow (X_{\zeta _j} , \mu _{\zeta _j})$ denote the factor map from $\alpha _{\zeta _k}$ to $\alpha _{\zeta _j}$, and let $\varphi _j : (X,\mu ) \rightarrow (X_{\zeta _j}, \mu _{\zeta _j})$ denote the factor map from $\alpha = \alpha _{\eta}$ to $\alpha _{\zeta _j}$. Applying the induction hypothesis to the action $\alpha _{\zeta _0}$ and the asymptotically equivariant sequence $\varphi _0 \circ \theta _n : (Y,\nu ) \rightarrow (X_{\zeta _0},\mu _{\zeta _0})$, we obtain a factor map $\theta ^0 :(Y,\nu )\rightarrow (X_{\zeta _0},\mu _{\zeta _0})$ and a subsequence $(n_i^{0})_{i\in\mathbb{N}}$ such that $\varphi _0 \circ \theta _{n_i^0}$ converges weakly to $\theta ^0$. Having defined the factor map $\theta ^j : (Y,\nu ) \rightarrow (X _{\zeta _j},\mu _{\zeta _j})$ and subsequence $(n_i^j)_{i\in \mathbb{N}}$, we apply the induction hypothesis to action $\alpha _{\eta _{j+1}}$ and the sequence $(\varphi _{j+1}\circ \theta _{n_i^j})_{i\in\mathbb{N}}$ to obtain a factor map $\theta ^{j+1} : (Y,\nu )\rightarrow (X_{\zeta _{j+1}},\mu _{\zeta _{j+1}})$ and a subsequence $(n_i^{j+1})_{i\in \mathbb{N}}$ of $(n_i^j)_{i\in \mathbb{N}}$ with $\varphi _{j+1} \circ \theta _{n_i^{j+1}}$ converging weakly to $\theta ^{j+1}$. Observe that $\varphi _{j,k}\circ \theta ^k = \theta ^j$ for all $j<k \in \mathbb{N}$. Since $\alpha =\varprojlim _j \alpha _{\zeta _j}$, this implies that there is a unique factor map $\theta =\varprojlim _j \theta ^j$ from $\beta$ to $\alpha$ such that $\varphi _j \circ \theta  = \theta ^j$ for all $j\in \mathbb{N}$. Morover, if we define the subsequence $(n_i)_{i\in\mathbb{N}}$ by $n_i = n^i_i$, then for each $j\in \mathbb{N}$, the sequence $( \varphi _j \circ \theta _{n_i} )_{i\in \mathbb{N}}$ converges weakly to $\theta ^j$, and hence the sequence $(\theta _{n_i})_{i\in\mathbb{N}}$ converges weakly to $\theta$. \hfill$\blacksquare$

{\it Proof of Theorem \ref{main}}.
This is immediate from Lemma \ref{lemma1} and Lemma \ref{lem:distal}.
\hfill$\blacksquare$

\begin{lemma}\label{lem:weakisom}
Let $\Gamma \curvearrowright ^{\alpha} (X,\mu )$ and $\Gamma \curvearrowright ^{\beta} (Y,\nu )$ be ergodic compact p.m.p.\ actions of $\Gamma$. If $\alpha$ are $\beta$ are weakly isomorphic, then they are isomorphic.
\end{lemma}

{\it Proof}. Let $\psi _0 :(X,\mu )\rightarrow (Y,\nu )$ and $\psi _1 :(Y,\nu )\rightarrow (X,\mu )$ be factor maps from $\alpha$ to $\beta$ and from $\beta$ to $\alpha$ respectively. It suffices to show that the map $\theta = \psi _1 \circ \psi _0$, factoring $\alpha$ onto itself, is an isomorphism. Since $\alpha$ is an ergodic compact action, we may assume that $(X,\mu ) = (K/L, m_{K/L})$, where $K$ is a compact metrizable group, $L<K$ is a closed subgroup, and that $\Gamma$ acts by translation on $K/L$ via a homomorphism $\Gamma \rightarrow K$ with dense image (in what follows we will identify $\Gamma$ with its image in $K$). In order to show that $\theta$ is injective on a conull subset of $K/L$, it is enough to show that the isometric linear embedding $T ^{\theta} : L^2(K/L ) \rightarrow L^2(K/L )$, $\xi \mapsto \xi \circ \theta$, which $\theta$ induces on $L^2(K/L )$, is surjective.

The translation action $K\curvearrowright K/L$ gives rise to a unitary representation $\lambda _{K/L}$, of $K$ on $\mathcal{H} = L^2 (K/L)$, which we may view as a subrepresentation of the left regular representation of $K$ of $L^2(K)$ via the inclusion $L^2(K/L) \hookrightarrow L^2(K)$ associated to the natural projection $K\rightarrow K/L$. The representation $\lambda _{K/L}$ may be expressed as a direct sum $\lambda _{K/L} = \bigoplus _{\pi \in \widehat{K}} \lambda _{K/L}^{\pi}$, where $\widehat{K}$ is a collection of representatives for isomorphism classes of irreducible unitary representations of $K$, and where for each $\pi \in \widehat{K}$, the representation $\lambda _{K/L}^{\pi}$ is the restriction of $\lambda _{K/L}$ to the closed linear span $\mathcal{H}^{\pi}$, of all subspaces of $\mathcal{H}$ on which $\lambda _{K/L}$ is isomorphic to $\pi$. For each $\pi\in\widehat{K}$ we may further write $\mathcal{H}^{\pi}$ as a direct sum of irreducible subspaces $\mathcal{H}^{\pi}=\bigoplus _{i<n_{\pi}}\mathcal{H}^{\pi , i}$, where for each $i<n_{\pi}$ the representation $\lambda _{K/L}^{\pi ,i} := \lambda _{K/L}|\mathcal{H}^{\pi ,i }$ is isomorphic to $\pi$. Since $\lambda _{K/L}$ is a subrepresentation of the left regular representation of $K$, by the Peter-Weyl Theorem, each $n_{\pi}$ is finite and therefore the subspaces $\mathcal{H}^{\pi}$ are all finite dimensional.

Since $\Gamma$ is dense in $K$, the operator $T^{\theta}$ intertwines $\lambda _{K/L}$ with itself, so by Schur's Lemma we see that $T^{\theta}$ intertwines each $\lambda _{K/L}^{\pi}$ with itself, i.e., $T ^{\theta} (\mathcal{H}^{\pi})\subseteq \mathcal{H}^{\pi}$ for all $\pi \in \widehat{K}$. Since $T^{\theta}$ is injective and each $\mathcal{H}^{\pi}$ is finite dimensional it follows that $T^{\theta} (\mathcal{H}^{\pi })= \mathcal{H}^{\pi}$ for all $\pi \in \widehat{K}$ and hence $T^{\theta}$ is surjective on $L^2(K/L )$, as was to be shown.
\hfill$\blacksquare$

{\it Proof of Corollary \ref{maincor1}}. Assume that $\alpha$ is weakly contained in $\beta$. Let $\alpha'$ be the maximal distal factor of $\alpha$. As $\alpha'$ is a factor of $\alpha$ and $\alpha \prec \beta$, we have that $\alpha' \prec \beta$. Since $\alpha'$ is distal and $\beta$ is strongly ergodic, by applying Theorem \ref{main} we deduce that $\alpha'$ is a factor of $\beta$, as claimed. \hfill$\blacksquare$

{\it Proof of Corollary \ref{maincor}}. Assume that $\alpha$ is weakly equivalent to $\beta$. Since $\beta$ is strongly ergodic, this implies that $\alpha$ is strongly ergodic. Let $\alpha '$ and $\beta '$ be the maximal distal factors of $\alpha$ and $\beta$, respectively.  By applying Corollary \ref{maincor1} twice we deduce that $\alpha'$ is a factor of $\beta'$ and that $\beta'$ is a factor of $\alpha'$. Thus, $\alpha '$ and $\beta '$ are weakly isomorphic. The moreover statement of Corollary \ref{maincor} now follows from Lemma \ref{lem:weakisom}.
\hfill$\blacksquare$

{\it Proof of Corollary \ref{cor3}}.  For simplicity, denote $\sigma=\alpha_{(a,b)}$ and $\sigma'=\alpha_{(a',b')}$.
Assume that $\sigma\prec\sigma'$. Denote by $A,B$ and $A',B'$ the matrices constructed from the pairs $(a,b)$ and $(a',b')$.
Since $A',B'$ have algebraic entries,  $\sigma'$ has spectral gap (see \cite[Theorem 1.2]{BdS14}). Corollary \ref{maincor} gives that $\sigma$ is a factor of $\sigma'$.
Let $\theta:K\rightarrow K$ be a measurable map such that $\theta(\sigma'(\gamma)x)=\sigma(\gamma)\theta(x)$, for all $\gamma\in\Gamma$ and almost every $x\in K$. If $t\in K$, then $K\ni x\rightarrow\theta(x)^{-1}\theta(xt)\in K$  is a $\sigma'(\Gamma)$-invariant map. Thus, there is $\delta(t)\in K$ such that $\theta(x)^{-1}\theta(xt)=\delta(t)$, for almost every $x\in K$. It follows that $\delta:K\rightarrow K$ is a continuous homomorphism and that there is $k\in K$ such that $\theta(x)=k\delta(x)$, for almost every $x\in K$. From this we deduce that $\delta(\sigma'(\gamma))=k^{-1}\sigma(\gamma)k$, for all $\gamma\in\Gamma$.

Since $K$ is a simple group and $\delta$ is non-trivial, $\delta$ must be one-to-one.  Since $K$ is not isomorphic to any of its proper closed subgroups, $\delta$ is also onto. Thus, since $K$ has no outer automorphisms, we can find $g\in K$ such that $\delta(x)=gxg^{-1}$, for all $x\in K$. Therefore, $g\sigma'(\gamma)g^{-1}=k^{-1}\sigma(\gamma)k$, for all $\gamma\in\Gamma$.
In particular, we get that $gA'a^{-1}=k^{-1}Ak$. Hence $A,A'$ must have the same trace, which implies that $\frac{a}{b}=\frac{a'}{b'}$.

Similarly, if $\alpha_{(a,b)}^n\prec\alpha_{(a',b')}^n$, for some $n\geqslant 2$, it follows that $\frac{a}{b}=\frac{a'}{b'}$. This completes the proof.
\hfill$\blacksquare$

{\it Proof of Corollary \ref{cor4}}. Assume by contradiction that there exists a p.m.p. action $\Gamma\curvearrowright^{\sigma}(Y,\nu)$ such that $\alpha\prec\sigma$ and $\sigma(\Gamma)(y)\subset\beta(\Lambda)(y)$, for almost every $y\in Y$.  Identify $\alpha$ with a left translation action $\Gamma\curvearrowright (K/L,m_{K/L})$ associated to a dense embedding of $\Gamma$ into a compact group $K$.
Let $d$ be a compatible metric on $X=K/L$. By Lemma \ref{lemma1} we can find a sequence $\theta_n:Y\rightarrow X$ of measurable maps such that $d(\theta_n(\gamma y),\gamma \theta_n(y))\rightarrow 0$, for all $\gamma \in\Gamma$ and almost every $y\in Y$. Denote by $\mathcal R_{\sigma}$ and $\mathcal R_{\beta}$ the equivalence relations associated to $\sigma$ and $\beta$, so that $\mathcal R_{\sigma}\subset\mathcal R_{\beta}$.

Since $\alpha$ is essentially free and $\alpha\prec\sigma$, $\sigma$ is essentially free. Since $\Gamma$ is non-amenable,  the restriction of $\mathcal R_{\sigma}$ to any non-negligible set $Y_0\subset Y$ is not hyperfinite. By \cite[Theorem 1]{CI08} we deduce the existence of a
$\sigma(\Gamma)$-invariant non-negligible measurable set $Y_1\subset Y$ such that the restriction $\sigma_1$ of $\sigma$ to $Y_1$ is strongly ergodic.
By applying Lemma \ref{equiv} to the restrictions of $\theta_n$ to $Y_1$ we conclude that $\alpha$ is a factor of $\sigma_1$. Let $\theta:Y_1\rightarrow X$ be a measurable, measure preserving map such that $\theta(\gamma y)=\gamma \theta(y)$, for every $\gamma \in\Gamma$ and almost every $y\in Y_1$.

We will reach a contradiction by applying \cite[Theorem 6.2]{CPS11}. To this end, we denote by $\Theta:L^{\infty}(X)\rtimes_{\alpha}\Gamma\rightarrow L^{\infty}(Y_1)\rtimes_{\sigma_1}\Gamma$ the $*$-homomorphism given by $\Theta(f)=f\circ\theta$ and $\Theta(u_{\gamma})=u_{\gamma}$, for all $f\in L^{\infty}(X)$ and $\gamma \in\Gamma$. We view $N:=L^{\infty}(Y)\rtimes_{\sigma}\Gamma$ as a subalgebra of $M:=L^{\infty}(Y)\rtimes_{\beta}\Lambda$. Let
$p={\bf 1}_{Y_1}\in N$, $P=\Theta(L^{\infty}(X))\subset pMp$ and $\mathcal G =\{au_{\gamma}p|a\in\mathcal U(L^{\infty}(Y)),\gamma \in\Gamma\}$.
Then $\mathcal G$ is a group of unitaries in $pMp$ which normalize $P$ and satisfies $\mathcal G''=Np$.

Next, denote by $\widetilde M=L(\mathbb Z\wr\Lambda)$ the von Neumann algebra of the wreath product group $\mathbb Z\wr\Lambda$.
Recall that $(Y,\nu)=(Z,\rho)^{\Lambda}$ and consider a fixed embedding of $L^{\infty}(Z)$ into $L(\mathbb Z)$. From this we get a $\Lambda$-equivariant embedding of  $L^{\infty}(Y)=L^{\infty}(Z)^{\Lambda}$ into $L(\oplus_{\lambda \in\Lambda}\mathbb Z)=L(\mathbb Z)^{\Lambda}$, and thus an embedding $M\subset\widetilde{M}$. It is easy to that $(L^{\infty}(Y)p)'\cap p\widetilde{M}p=L(\oplus_{\lambda \in\Lambda}\mathbb Z)p$.

Thus, if we denote $Q=\mathcal G'\cap p\widetilde Mp$, then $Q\subset L(\oplus_{\lambda \in\Lambda}\mathbb Z)p$. Since $Q$ commutes with $Np$ and $N$ has no amenable direct summand, \cite[Theorem 2]{CI08} implies that $Q$ is completely atomic.

Let $q\in Q$ be a non-zero projection such that $Qq=\mathbb Cq$.
Then $(\mathcal Gq)'\cap q\widetilde Mq=\mathbb Cq$.
Morever, since $\alpha$ is a compact action, the conjugation action $\mathcal G\curvearrowright P$ is compact. Hence, the conjugation action $\mathcal Gq\curvearrowright Pq$ is compact and thus weakly compact (see \cite[Definition 6.1]{CPS11}). By applying \cite[Theorem 6.2]{CPS11} we get that either $(\mathcal Gq)''=Nq$ has an amenable direct summand or $P\prec_{\widetilde M}L(\Lambda)$ (see \cite[Definition 2.1]{CPS11}).
Since $N$ has no amenable direct summand and $P\subset L(\oplus_{\lambda \in\Lambda}\mathbb Z)$, neither of these conditions hold true, which gives the desired contradiction.
\hfill$\blacksquare$

\end{document}